\documentclass[10pt]{amsart}
\usepackage[utf8]{inputenc}
\usepackage[english,russian]{babel}

\usepackage{amsmath}
\usepackage{amssymb}
\usepackage{amsfonts}

\newtheorem{lemma}{Lemma}
\newtheorem{theorem}{Theorem}

\begin{document}
\renewcommand{\refname}{References}
\begin{center}
\textbf {ON SOME NEW EMBEDDINGS IN MINIMAL BOUNDED HOMOGENEOUS DOMAINS IN $\mathbb C^n$} \\

\vspace{8pt}

\textbf{R.F. Shamoyan, N.M. Makhina}\\
\end{center}

\vspace{8pt}

\noindent{\small {\bf Abstract.}} In this short note we consider very general bounded minimal homogeneous domains. Under certain natural additional conditions new sharp results on Bergman type analytic spaces in minimal bounded homogeneous domains are obtained. Domains we consider here are direct generalizations of the well-studied so-called bounded symmetric domains in $\mathbb C^n.$ In the unit disk and in the unit ball  all our results were obtained by first author. Some results were obtained previously in tubular and bounded pseudoconvex domains. Our proofs are heavily based on properties of so called  r-lattices for these general domains provided in recent papers of Yamaji. Our proofs are also based on   arguments provided earlier in less general domains.We partially extend an embedding result from Yamaji's paper.

\vspace{10pt}


\section{\textbf{Introduction}}
The goal of this paper to obtain new sharp results on Bergman type analytic
spaces in minimal bounded homogeneous domains. Our results previously were
known only in very particular case of such type domains in the unit ball. Some
results are known in tube domains and pseudoconvex domains (see~\cite{2}, \cite{4}-\cite{6}).
Tube domains are unbounded, pseudoconvex domains are not symmetric. Our
results are heavily based on a series of subtle new estimates obtained recently in~\cite{9}-\cite{12}. We note domains we consider here are direct generalizations of the
well-studied so-called bounded symmetric domains in $\mathbb C^n$ (see~\cite{9}-\cite{12}). Note, also, all mentioned domains and even polydisk are examples of minimal domains.
To formulate our results we need some basic notations. We say the bounded $\mathcal U$ domain in $\mathbb C^n$ is a minimal domain with a center $t \in \mathcal U$ if the following condition
is satisfied for every biholomorphism $\psi:\mathcal U \in \mathcal {U'}$ with $det J(\psi, t) = 1$ we have $Vol(\mathcal (U') \ge Vol(\mathcal U)$ where $J(\psi, t)$ denotes the complex Jacobi matrix of $\psi$ at $t$ (see \cite{9}-\cite{12}). We denote constants by $C, C_1, . . .$ We fix a minimal bounded homogeneous domain $\mathcal U$ with center $t$. 

Let $dV_\beta(z) = K_{\mathcal U}(z, z)^{-\beta}dV(z),$ $\beta \in \mathbb R,$ and $dV$ denote the Lebesgue measure on $\mathcal U$ (see \cite{9}-\cite{12}). Let $L^p_{\alpha, \beta}(\mathcal U, dV_\beta) = L^p(\mathcal U, dV_\beta) \cap H(\mathcal U),$ $0 < p \le \infty,$ where $H(\mathcal U)$ is a class of all analytic functions on $\mathcal U$. The last spaces are non-trivial if and only if $\beta > \beta_{min}$ for some fixed $\beta_{min}$ (see \cite{9}-\cite{12}). Note we will always assume this bellow. These are Banach spaces for $p\ge 1.$
We bellow denote by $K^{(\beta)}_{\mathcal U}$ the reproducing kernel of $L_\alpha^2(\mathcal U, dV_\beta).$ $L^2_\alpha$ is the Bergman spaces on $\mathcal U$  (unweighted) and $L_\alpha^2(\mathcal U, dV) = L^2(\mathcal U, dV)\cap H(\mathcal U)$. It is known that $K_\beta = K^{(\beta)}_\mathcal U (z, \omega) = C_\beta K_{\mathcal U}(z, \omega)^{1+\beta}$ for some positive constant $C_\beta$ (see ~\cite{9}-\cite{12}).

Yamaji (see \cite{9}) give criteria for the boundedness of positive Toeplitz
operators on weighted Bergman spaces of a minimal bounded homogeneous
domain in terms of the Berezin symbol or the averaging function of the symbol. Moreover, Yamaji estimate the essential norm of positive Toeplitz operators assuming that they are bounded. As an application of these estimates, also give necessary and sufficient conditions for the positive Toeplitz operators to be compact.

Using an integral formula on a homogeneous Siegel domain, Yamaji (see \cite{10}) give a necessary and sufficient condition for composition operators on the weighted Bergman space of a minimal bounded homogeneous domain to be compact in terms of a boundary behavior of the Bergman kernel.

In his paper \cite{10} Yamaji intrоduced the so called Berezin symbol and the averaging function for a positive Borel measure $\mu$ on a minimal bounded homogeneous domain $\mathcal U.$ These concepts were used by him in \cite{10} to get a complete description of all so called  Carleson measures
for $L^p_a(\mathcal U, dV_{\beta})$ analytic Bergman type spaces in minimal homogeneous domains for all positive values of $p.$ Namely it was shown
in \cite{10} that both Berezin symbol and  the averaging function in $\mathcal U$ must be bounded in $\mathcal U$ if and only if the positive Borel measure $\mu$ is a Carleson measure in $\mathcal U$ (see definition of Carleson measure in \cite{10}). For analytic Bergman type function spaces we mentioned above. 

We continue in this short note his investigation using  and modifying his arguments providing various new embeddings in minimal bounded homogeneous domains for various new mixed norm Bergman spaces.

\section{\textbf{Main results}}

The Bergman kernel $K_{\mathcal U}(z, \omega)$ of $\mathcal U$ is playing very important in our theorems below. Let $d_\mathcal U (\cdot, \cdot)$ be the Bergman distance on $\mathcal U.$ For any $z \in \mathcal U,$ $r > 0,$ let also $B(z, r) = \{\omega \in \mathcal U: d_{\mathcal U}(z, \omega) \le r\}$ be the Bergman metric disk with center $z$ and radius $r.$

The following result is a sharp result of distances on such domains (see~\cite{4}-\cite{6} for other domains). The proof is close to the case of unit disk, based on estimates from~\cite{9}-\cite{12} for such domains.

We formulate below  our all main results of this short note namely new embedding theorems for minimal bounded homogeneous domains for some new mixed norm analytic Bergman type spaces. The first result contains only sufficient condition on positive Borel measure, the other two results are sharp embeddings. Proofs use arguments provided earlier in less general domains and properties of so called sampling sequences(see lemmas below). In all theorems we assume that Forelly-Rudin estimates are valid (see lemma 4 and the remark after it.) 

The existence of so-called Bergman sampling sequence can be seen in~\cite{9}-\cite{12} (see also Lemma 3 below). This sequence and estimates of Bergman kernel on $\{B(z_k, \rho)\}$ balls are very vital for this paper. We denote below the Lebesques measure of $B(z, \rho)$ ball by $Vol.$ We denote by $Vol(E)$ the volume of $E$ set.

\begin{theorem}
Let $1\le q,s,r<\infty,$ $q\ge s,$ $\alpha>-1.$ Let $\{z_k\}$ be a sampling sequence in $\mathcal U,$ $\mu$ be a positive Borel measure in $\mathcal U.$ Then 
\[\int_{\mathcal U}|f(z)|^q d\mu(z)\le C\int_{\mathcal U}\left(\int_{B(\omega,r)}|f(z)|^sdV_{\alpha}(z)\right)^{q/s}dV(\omega),\]
if $\mu(B(z_k,r))\le C(Vol(B(z_k,r))^{\tau},$ for all $\{z_k\},$ $k=1,...,$ and for some positive $\tau,$ $\tau=\tau(s,\alpha,q).$
\end{theorem}

The following theorem is a new sharp result on embeddings in $L^2_{a}
(\mathcal U, dV_\beta)$ type analytic function spaces.

\begin{theorem}
Let $1 < p,q < \infty,$ $ 1 < s \le p < \infty,$ $\beta > \beta_{min},$ $ \rho > 0.$ Let $\{z_k\}$ be a sampling sequence in $\mathcal U,$ let $\mu$ be positive Borel measure on $\mathcal U$. Then
\[\sum _{k=1}^{\infty }\left(\int _{B(z_{k} ,\rho )}\left|f(z)\right|^{q} d\mu (z) \right)^{\frac{p}{q} } \le C\left\| f\right\| _{L_{\alpha }^{s} (\mathcal U,dV_{\beta } )} \eqno(1)\] 
if and only if
\[\mu(B(z_k,\rho))\le C(Vol)(B(z_k,\rho))^{\beta_0},\eqno(2) \]
for all $\{z_k\}\in \mathcal U,$ $\rho>0,$ $k=1,2,3,...$ for some fixed $\beta_0,$ $\beta_0=\beta_0(p,q,s,\beta).$
\end{theorem}

Note that the sufficiency part in this theorem follows immediately from lemmas 2-3 which we formulated below as in case ofmuch simpler domain   the unit ball (see~\cite{2}, \cite{8}).

Note it was shown in~\cite{9}-\cite{12} that the similar to (2) condition holds if and only if
\[\int _{\mathcal U}\left|f(z)\right|^{p} d\mu (z) \le \tilde{C}\int _{\mathcal U}\left|f(z)\right|^{p} dV_{\beta } (z) ,\eqno(3)\] 
for all $p > 0$ and for all $f \in L^p_{a}(\mathcal U, dV_{\beta}).$

Note the proofs of theorems of this paper can be obtained after careful study of estimates of the proof of the unit ball case and parallel estimates obtained recently in the case of bounded minimal homogeneous domains in $\mathbb C^n$ (see~\cite{2}, \cite{4}-\cite{6}, \cite{9}-\cite{12}).

The following theorem is another new sharp result on embeddings in $L^p_{a}(\mathcal U, dV_{\beta})$ analytic function spaces in minimal bounded homogeneous domain in $\mathbb C^n.$ The base of proof is Forelly - Rudin estimates and lower estimate for Bergman kernel (which we assume are valid for our general domains).

\begin{theorem}
Let $\mu$ be a positive Borel measure on $\mathcal U,$ and $\{z_k\}$ be a Bergman sampling sequence. Let $\alpha > \alpha_{min},$ $f_i \in H(\mathcal U),$ $1 < p_i, q_i < \infty,$ $i = 1, . . . , m,$ so that
$\sum\limits_{i=1}^m {\frac{1}{q_i}}=1.$ Then
\[\int _{\mathcal U}\prod _{i=1}^{m}\left|f_{i} (z)\right|^{p_{i} }  d\mu (z) \le \tilde{\tilde{C}}\prod _{i=1}^{m}\left[\sum _{k=1}^{\infty }\left. \left(\int _{B(z_{k} ,2r)}\left|f_{i} (z)\right|^{p_{i} } dV_{\alpha } (z) \right)^{q_{i} } \right] ^{\frac{1}{q_{i} } } \right. \eqno(4) \] 
if and only if $\mu(B(z_k, r)) \le C(Vol(B(z_k, r)))^{\alpha_0}$ for every $k, k = 1, 2, . . . ,$ $r > 0,$ for some fixed $\alpha_0, \alpha_0(p_1..,p_m, q_1,...q_m,\alpha,m).$
\end{theorem}

Note that the sufficiency part in this theorem follows immediately from lemmas 2-3 which we formulated below as in case ofmuch simpler domain   the unit ball (see~\cite{2}, \cite{8}).

Very similar results with similar proofs were obtained by the first author in
tubular domains over symmetric cones (unbounded domains) and bounded strictly
pseudoconvex (nonsymmetric) domains (see~\cite{1}, \cite{2}, \cite{8}, and references there).

Some words on proofs:

We supply three lemmas from~\cite{9}-\cite{12} which are crucial for proofs. Analogues in tube and pseudoconvex domains can bee seen in~\cite{1}, \cite{2}, \cite{8}.

Below in lemma 1 we provide an important for our proofs of theorem 2, and 3 estimate from below of Bergman kernel on Bergman balls (see~\cite{8}-\cite{12}). Note also that they are also valid for many other general domains in $\mathbb C^n.$

\begin{lemma}
(see~\cite{12}) Take $\rho > 0.$ Then there exists $C_{\rho} > 0$ such that
\[C_{\rho }^{-1} \le \left|\frac{K_{\mathcal U} (z,a)}{K_{\mathcal U} (a,a)} \right|\le C_{\rho } ,z,a\in \mathcal U,\beta _{\mathcal U} (z,a)\le \rho ,\eqno(5)\] 
where $\beta_{\mathcal U}$ means the Bergman distance on $\mathcal U.$
\end{lemma}

Note that these estimates are valid also for weighted Bergman kernel (see~ \cite{9}-\cite{2}).

\begin{lemma}
(see~\cite{9}) There exists a positive constant $C$ such that
\[\left|f(a)\right|^{p} \le \frac{C}{Vol(B(a,\rho ))} \int _{B(a,\rho )}\left|f(z)\right|^{p} dV(z), \eqno(6)\] 
$f \in H(\mathcal U),$ $ p > 0,$ $ a \in \mathcal U.$
\end{lemma}

This lemma is valid also in more general form when we replace $dV$ by $dV_\beta$ (see~\cite{8}-\cite{12}).

\begin{lemma}
(see~\cite{12}) There exists a sequence $\{\omega_j\}\in \mathcal U$  satisfying the following conditions
\[\mathcal U=\bigcup _{j=1}^{\infty }B(\omega _{j} ,\rho ),B(\omega _{i} ,\rho /4)\cap B(\omega _{j} ,\rho/4 ) =\emptyset ,i\ne j.\] 
There exists a positive integer $n$ such that for each point $z \in \mathcal U$ belongs to at most $n$ of sets $B(\omega_j, 2\rho).$
\end{lemma}

\begin{lemma}
(Forelly-Rudin estimates in $\mathcal U$ (see~\cite{8}-\cite{12}))
We have 
\[\int_{\mathcal U}|K_{\mathcal U}(z,z'')|^{1+\alpha}dV_\beta(z'')< c(\alpha,\beta)|K_{\mathcal U}(z,z)|^{\alpha-\beta},\]
$\alpha=1+2\beta,$ $\beta>\beta_0,$ $z \in \mathcal U,$ $\beta_0$ is large enough.
\end{lemma}

We assume that this lemma is valid in more general situation when indexes $\alpha$ and $\beta$ are independent.

All results of this note with very similar proofs in context of bounded strongly pseudoconvex domains with smooth boundary were proved earlier by the first author.

We obtained similar results concerning so-called distance function $dist_f(g, X)$ in these bounded minimal domains and they will be presented by us in another paper. Similar extremal (distance) problems and results for other domains were proved by the first author in \cite{1}, \cite{3}-\cite{6}.

\end{document}